\documentclass[11pt,reqno,preprint]{elsarticle}

\makeatletter
\def\ps@pprintTitle{%
 \let\@oddhead\@empty
 \let\@evenhead\@empty
 \def\@oddfoot{}%
 \let\@evenfoot\@oddfoot}
\makeatother

\usepackage{amsmath,amsthm,amscd,wrapfig,amsfonts,amssymb,mathtools}
\usepackage{enumitem,comment,mathabx,mathrsfs}
\usepackage{booktabs,multirow,lscape,cool}
\usepackage{url,parskip,caption}
\usepackage[top=1.1in, bottom=1.1in, left=1.1in, right=1.1in]{geometry}

\newtheorem{remark}{Remark}

\providecommand{\PS}[1]{\ensuremath{\mathit{Ps}\hskip.05em{}_{#1}}}
\providecommand{\QS}[1]{\ensuremath{\mathit{Qs}\hskip.05em{}_{#1}}}

\providecommand{\ps}[1]{\ensuremath{\mathsf{Ps}\hskip.05em{}_{#1}}}
\providecommand{\qs}[1]{\ensuremath{\mathsf{Qs}\hskip.05em{}_{#1}}}

\providecommand{\p}[1]{\ensuremath{\mathsf{P}\hskip.05em{}_{#1}}}
\providecommand{\q}[1]{\ensuremath{\mathsf{Q}\hskip.05em{}_{#1}}}

\providecommand{\PsiS}[1]{\ensuremath{\mathit{\Psi S}\hskip.05em{}_{#1}}}

\providecommand{\e}[1]{\ensuremath{\times 10^{#1}}}

\makeatletter
\g@addto@macro\normalsize{%
  \setlength\abovedisplayskip{.4em}
  \setlength\belowdisplayskip{.4em}
  \setlength\abovedisplayshortskip{.4em}
  \setlength\belowdisplayshortskip{.4em}
}

\begin{document}

\begin{frontmatter}

\begin{abstract}
The standard algorithm for the numerical evaluation of
 the prolate spheroidal wave function $\ps{n}(x;\gamma^2)$
of order $0$, bandlimit $\gamma > 0$  and characteristic exponent $n$
has  running time which grows with both $n$ and $\gamma$.
Here, we describe an alternate approach
which runs in time independent of these quantities.
We present the results of numerical experiments demonstrating the properties
of our scheme, and we have made our implementation of  it publicly available.



\end{abstract}


\title
{
An $\mathcal{O}(1)$ algorithm for the numerical evaluation of
the prolate spheroidal wave functions of order $0$
}

\author[1]{Xinge Zhang}

\author[1]{James Bremer\corref{cor1}}
\ead{bremer@math.ucdavis.edu}

\cortext[cor1]{Corresponding author}

\address[1]{Department of Mathematics, University of California, Davis}








\end{frontmatter}

The  prolate spheroidal wave functions
\begin{equation}
\ps{0}\left(z;\gamma^2\right),
\ps{1}\left(z;\gamma^2\right),
\ps{2}\left(z;\gamma^2\right),
\ldots
\label{introduction:ps}
\end{equation}
are the eigenfunctions of the restricted Fourier operator
\begin{equation}
T_\gamma\left[f\right](z) = \int_{-1}^1 \exp(i \gamma z t) f(t) \ dt.
\label{introduction:intop}
\end{equation}
As such, they provide an efficient mechanism for representing bandlimited functions,
and for performing many computations related to such functions
(see, for instance, 
\cite{PSWFI, PSWFII, PSWFIII, Osipov-Rokhlin-Xiao,Hogan-Lakey,SHKOLNISKY2007235,Rhodes,REYNOLDS2013352,beylkin1}).

In the seminal work \cite{PSWFI}, it was observed
these functions also constitute the set of solutions of a singular 
self-adjoint Sturm-Liouville problem.  More explicitly, (\ref{introduction:ps}) is 
the collection of all eigenfunctions of the  differential operator
\begin{equation}
L_\gamma\left[y\right](x) = - (1-x^2)y'(x) + 2 x y'(x) + \gamma^2 x^2 y(x)
\label{introduction:diffop}
\end{equation}
which satisfy the self-adjoint boundary conditions
\begin{equation}
\lim_{x\to 1} (1-x^2) y'(x)  = 0 = \lim_{x\to -1} (1-x^2) y'(x).
\label{introduction:bc}
\end{equation}
In other words, each of the functions $\ps{n}(x;\gamma^2)$ is a solution
of the spheroidal wave equation
\begin{equation}
(1-x^2)y'(x) - 2 x y'(x) + (\chi -  \gamma^2 x^2) y(x) = 0
\label{introduction:swe}
\end{equation}
with  $\chi$  the Sturm-Liouville eigenvalue of $\ps{n}(x;\gamma^2)$;
we will denote this value of $\chi$ via  $\chi_n(\gamma^2)$.

The standard algorithm for the numerical evaluation
of the functions (\ref{introduction:ps})
is based on discretizing the eigenproblem for the differential operator $L_\gamma$  rather than the
eigenproblem for the integral operator $T_\gamma$.    
This approach is preferred because of the nature of the spectrum of  $T_\gamma$  
 --- it  has roughly $2\gamma/\pi$ 
eigenvalues which are close to  $\sqrt{2\pi/\gamma}$, on the order of $\log(\gamma)$ eigenvalues
which decay rapidly to $0$ and the rest of its 
eigenvalues are of very small magnitude.  The eigenvalues  of $L_\gamma$
are, on the other hand, well-separated.

When the eigenproblem 
\begin{equation}
L_\gamma\left[y\right](x) = \chi y(x)
\label{introduction:eigenproblem}
\end{equation}
is discretized by introducing  one of the representations
%
%
\begin{equation}
y(x) = \sum_{k=0}^\infty a_{k} \p{2k}(x)
\label{introduction:legendre_expansion1}
\end{equation}
or
\begin{equation}
y(x) = \sum_{k=0}^\infty a_{k} \p{2k+1}(x),
\label{introduction:legendre_expansion2}
\end{equation}
where $\p{n}(x)$ denotes Ferrer's version of the Legendre function
of the first  kind of degree $n$,
the result is an infinite symmetrizable tridiagonal  matrix.
In the case of (\ref{introduction:legendre_expansion1}),
the eigenvalues of the matrix are 
%
\begin{equation}
\chi_0(\gamma^2) < \chi_2(\gamma^2) <  \chi_4(\gamma^2) < \cdots
\label{introduction:pseven}
\end{equation}
while the representation (\ref{introduction:legendre_expansion2}) 
leads to a matrix whose spectrum consists of 
%
\begin{equation}
\chi_1(\gamma^2) <  \chi_3(\gamma^2)<  \chi_5(\gamma^2) < \cdots .
\label{introduction:psodd}
\end{equation}
After proper normalization, the eigenvectors of these matrices give the coefficients
in the Legendre expansions of the functions (\ref{introduction:ps}).
To evaluate $\ps{n}(x;\gamma^2)$,  its Legendre expansion is first constructed
by truncating  one of the two infinite matrices mentioned above, symmetrizing it,
and  calculating the appropriate eigenvalue and corresponding eigenvector.
The fact that the matrix is symmetrizable greatly reduces the difficulty of these computations.
The function $\ps{n}(x;\gamma^2)$ can then be evaluated at any point
in the interval $(-1,1)$ via the resulting Legendre expansion.

This procedure is sometimes called the Legendre-Galerkin method,
although we  refer to it as the Xiao-Rokhlin algorithm
as it appears to have first been described in its entirety in
\cite{Xiao-Rokhlin-Yarvin}.
A part of the procedure was described earlier in \cite{Hodge}; 
there, the Sturm-Liouville eigenvalues
of the prolate spheroidal wave functions are obtained in the fashion described above.
However, the three-term recurrence relations satisfied by the Legendre coefficients
are  then used to construct the expansions
of the prolate spheroidal wave functions.   
The resulting method does not provide a numerically stable mechanism for evaluating
the functions $\ps{n}(x;\gamma^2)$  and
extended precision arithmetic is required
to produce accurate results using it, even for small values
of $\gamma$ and $n$.
It is often stated in the literature
that the Xiao-Rokhlin procedure was first described in  \cite{Bouwkamp}.
In fact, the method of \cite{Bouwkamp} for the computation
of the Sturm-Liouville eigenvalues is quite different and
is based on the well-known observation that the three-term recurrence relations
which arise from inserting the representations (\ref{introduction:legendre_expansion1})
and  (\ref{introduction:legendre_expansion2}) into (\ref{introduction:swe}) 
are related to infinite continued fraction expansions.
A thorough discussion of the Xiao-Rokhlin algorithm and related techniques
can be found in \cite{Osipov-Rokhlin-Xiao}.

We refer to the construction of the Legendre expansion as the precomputation phase
of the Xiao-Rokhlin algorithm.   Its running time clearly depends on the 
size of the  tridiagonal discretization matrix formed during the procedure.
The  precise dimension  required to achieve a specified precision for 
$\ps{n}(x;\gamma)$ is not known; however,
the numerical experiments described in \cite{Feichtinger}
suggest that it must grow as 
\begin{equation}
\mathcal{O}\left( n + \sqrt{n \gamma} \right)
\label{introduction:cost}
\end{equation}
in order to achieve fixed precision independent of $\gamma$ and $n$.
It is likely, then,  that the cost of the Xiao-Rokhlin precomputation phase
is 
\begin{equation}
\mathcal{O}\left( \left(n + \sqrt{n \gamma}\right)\log\left(n + \sqrt{n \gamma}\right) \right)
\end{equation}
since  one eigenvalue and eigenvector pair of an $n \times n$ symmetric
tridiagonal matrix can be found in $\mathcal{O}(n \log (n))$ operations.
Moreover, the cost of evaluating the resulting Legendre expansion
at a single point scales linearly with the dimension of the tridiagonal
matrix so that the asymptotic complexity of evaluating $\ps{n}(x;\gamma^2)$
at a point via the Xiao-Rokhlin method 
is likely $\mathcal{O}\left( n + \sqrt{n \gamma} \right)$, exclusive of the costs
of the precomputation phase.
Of course, there exist fast algorithms for evaluating an $n$-term Legendre expansion  
at certain collections of $\mathcal{O}(m)$ points in $\mathcal{O}\left(n \log(n) + m\right)$ operations.
However, it is often desirable to evaluate the prolate functions at a small
number of points.   This is particularly true when performing calculations
on parallel computers, in which case it is often preferable to perform
separate evaluations  on different computational units.  
Hence, the cost of evaluating $\ps{n}(x;\gamma^2)$ at a 
single point is of interest.

Here, we describe an algorithm for the numerical evaluation of 
the prolate spheroidal wave functions
  that runs in time independent  of $n$ and $\gamma$.  
Like the Xiao-Rokhlin algorithm, it has a precomputation phase
in which a representation of  $\ps{n}(x;\gamma^2)$
is constructed.  However, rather than a Legendre expansion, we use
a  nonoscillatory phase function  $\PsiS{\nu}(x;\gamma^2)$ for the differential
equation (\ref{introduction:swe}) to represent $\ps{n}(x;\gamma^2)$.
 The cost of constructing $\PsiS{\nu}(x;\gamma^2)$
is independent of the  parameters $n$ and $\gamma$, as is the cost
of evaluating $\ps{n}(x;\gamma^2)$ via $\PsiS{\nu}(x;\gamma^2)$.

That such nonoscillatory phase functions exist for  many second order differential
equations, including (\ref{introduction:swe}), has long been known,
and they have often been exploited to accelerate
the numerical calculation of certain special functions  and  their zeros.
For instance,   \cite{olver_1950} introduces a scheme
for  calculating the zeros of Bessel functions which also relies on 
the existence of a nonoscillatory phase function for Bessel's differential equation,
\cite{Goldstein-Thaler} suggests the use of such a phase function 
 to evaluate the Bessel functions
of large orders,  and 
one component of the widely used algorithm of \cite{Amos} for the evaluation of Bessel
functions makes use of asymptotic expansions of a nonoscillatory phase function
for Bessel's differential equation.   
These algorithms rely on extensive analytic understanding of the nonoscillatory
phase functions for Bessel's differential equation.
Similar approaches can only be applied to other classes of special function
for which such extensive knowledge is available, and this has, until recently,
limited the applicability of such methods to a relatively  small number of cases.

In \cite{SpiglerZeros}, an iterative method based on nonoscillatory phase functions is used to
calculate the zeros of functions
defined by a class of second order differential equations with polynomial coefficients.
A generalization of the scheme is used in \cite{SpiglerPhase1} to construct
asymptotic expansions for the 
 solutions of a large class of second order differential equations
with polynomials coefficients, and \cite{SpiglerPhase2} extends the method
to second order differential equations with coefficients that grow
exponentially fast.  The schemes of \cite{SpiglerZeros}, \cite{SpiglerPhase1}  
and \cite{SpiglerPhase2} are widely applicable, but they require the calculation
of high order derivatives of the coefficients of the differential
equations to which they are applied.  Since these derivatives cannot be
calculated numerically without severe loss of precision, the schemes of
\cite{SpiglerZeros}, \cite{SpiglerPhase1}   and \cite{SpiglerPhase2}
are carried out symbolically using computer algebra systems.
They are best viewed as algorithms for the symbolic computation of 
asymptotic expansions of nonoscillatory phase functions.  
 
In  \cite{BremerKummer}, a numerical algorithm for the calculation of 
nonoscillatory phase functions for a large class of second order
differential equations is introduced.  It runs in time independent 
of the frequency of oscillation of their solutions 
and does not require knowledge of the derivatives of the coefficients
of the differential equation.    This algorithm is a major component
in our  scheme to construct the nonoscillatory phase function
for (\ref{introduction:swe}).

By itself, however, the algorithm of \cite{BremerKummer} is insufficient to 
construct $\PsiS{n}(x;\gamma^2)$   since it requires knowledge of the
 Sturm-Liouville eigenvalue $\chi_n(\gamma^2)$ to do so.
The obvious solution is to use an asymptotic expansion
for $\chi_n(\gamma^2)$.
 Unfortunately, high accuracy  asymptotic expansions of the Sturm-Liouville eigenvalues
which are suitable for use in  numerical codes  
do not appear to be available at this time.  
For instance, \cite{Dunster} and \cite{Bonami2016} 
describe   Liouville-Green type uniform asymptotic
expansions of $\chi_n(\gamma^2)$ and $\ps{n}(x;\gamma^2)$, but they involve
a complicated change of variables defined in terms of an
elliptic integral.  Moreover,  only a few coefficients in the resulting 
expansions are known, which greatly limits the achievable
accuracy of these expansions.

The construction of asymptotic expansions for $\chi$ and 
the prolate spheroidal wave functions is complicated by the fact that,
when viewed as a function of the characteristic exponent $\nu$,
$\chi$  has branch cuts.
We give a definition of the  characteristic exponent of a solution of (\ref{introduction:swe})
in Section~\ref{section:charexp};   for the reader who is unfamiliar
with the notion, it suffices for now to know that it 
describes the behavior of the solution at infinity and allows for the 
extension of the definitions of  $\ps{\nu}(x;\gamma^2)$,
$\PsiS{\nu}(x;\gamma^2)$ and $\chi_{\nu}(\gamma^2)$
to arbitrary complex values of $\nu$.  Many asymptotic expressions for
$\chi$ and other quantities related to the prolate spheroidal wave functions
are, in effect, expansions in $\nu$, and so are attempts to
 approximate discontinuous functions with all of the obvious difficulties 
that entails.

\begin{figure}[t]
\begin{center}
\includegraphics[width=.49\textwidth]{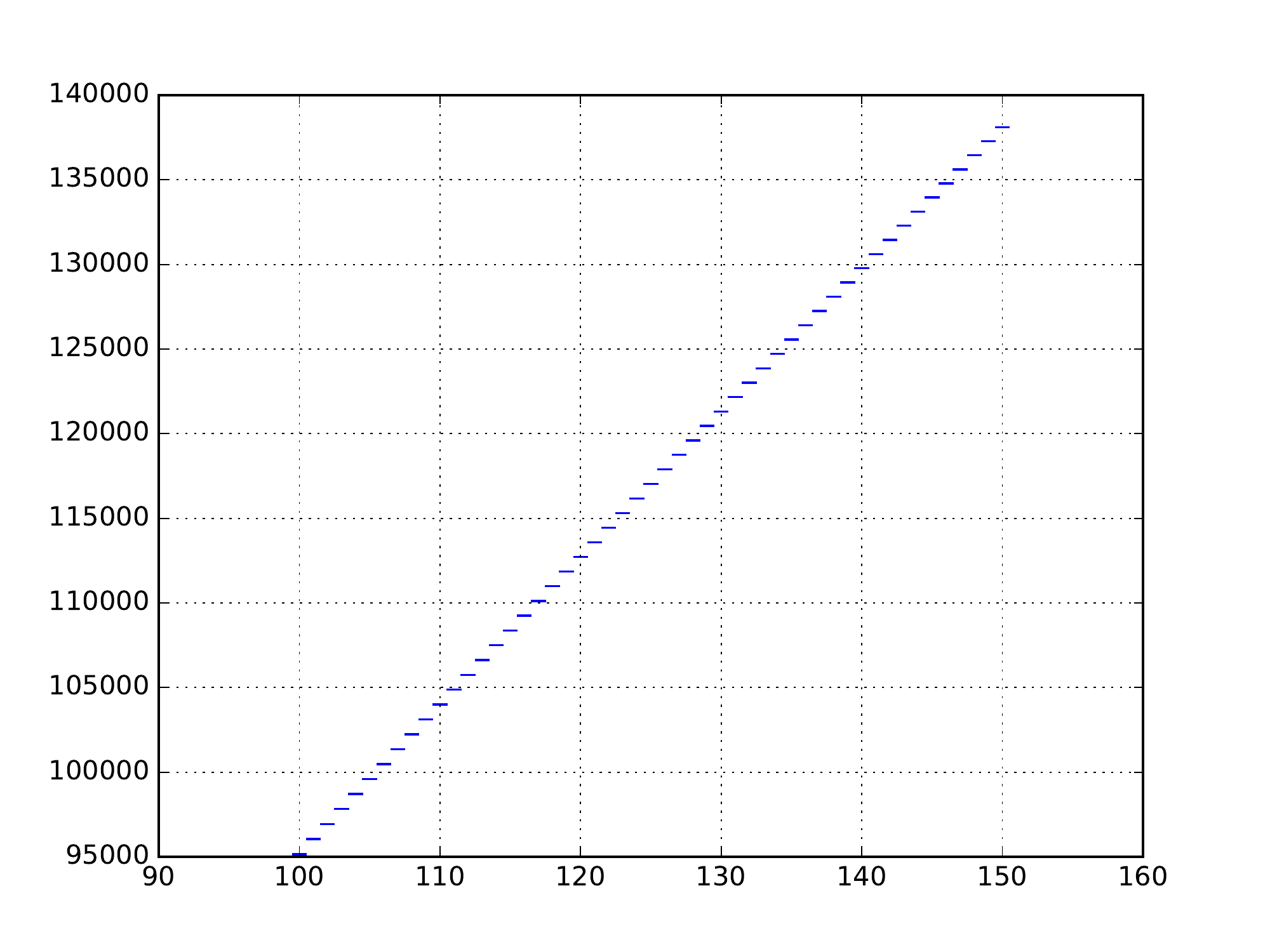}
\hfil
\includegraphics[width=.49\textwidth]{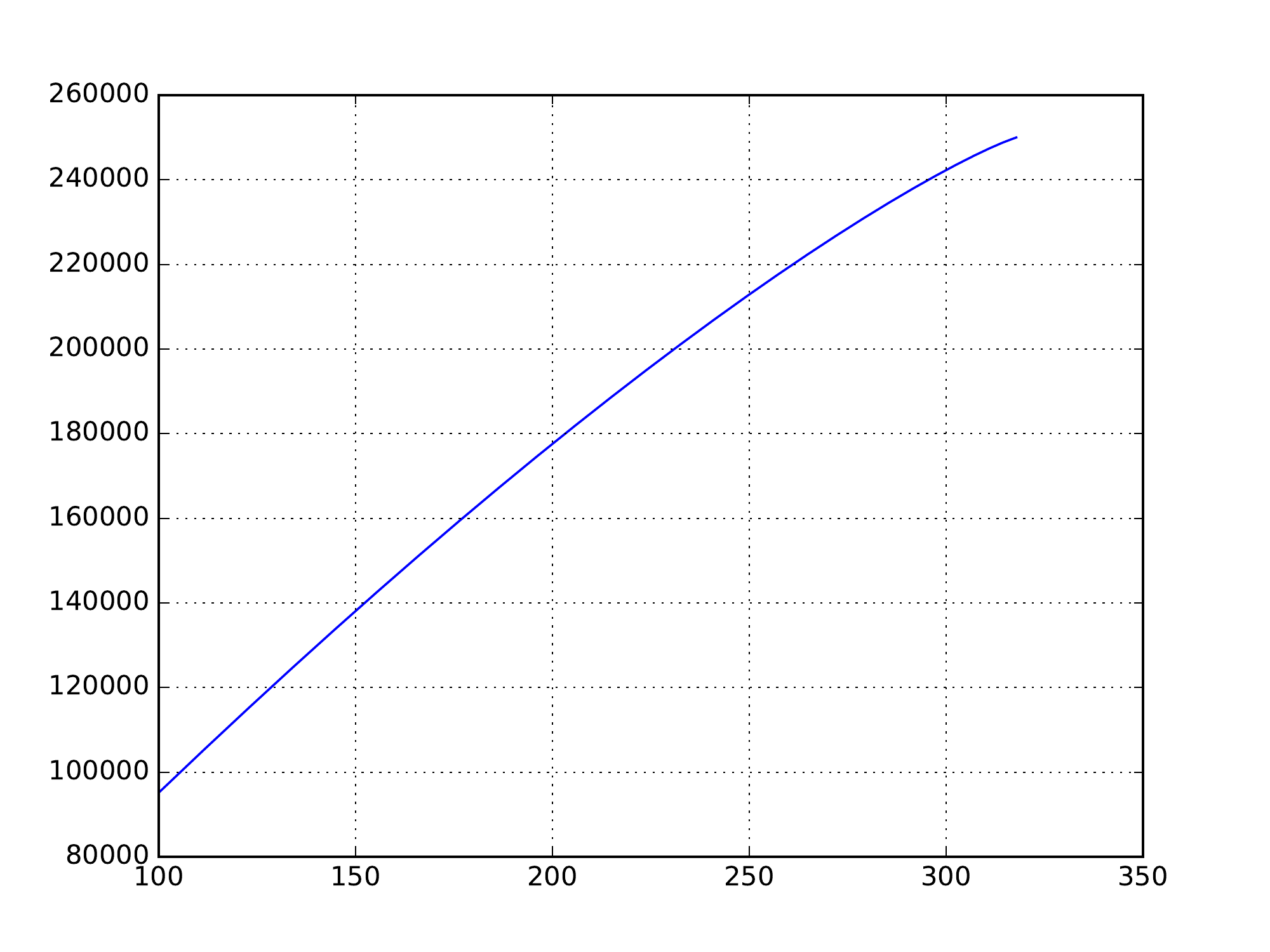}
\end{center}
\caption{On the left is a plot of $\chi$ as a function of the characteristic exponent
$\nu$ with $\gamma$ fixed at $500$.
There is a discontinuity in the graph of this function at each half-integer value
of $\nu$.
On the right is a plot of $\chi$ as a function of $\xi$
with  $\gamma$ again fixed at $500$.
The parameter $\xi$
to the value $\PsiS{n}(0;\gamma^2)$ of a nonoscillatory phase function
for the spheroidal  wave equation at $0$ through an affine mapping. }
\label{figure:chiplots}
\end{figure}

We take a somewhat different tack and numerically construct expansions
of $\chi$ and and a few other quantities  as functions of the bandlimit $\gamma$ and
a second parameter which is related to $\PsiS{n}(0;\gamma^2)$, the value at $0$ of a nonoscillatory phase
function representing $\ps{n}(x;\gamma^2)$.  It is the case
that
\begin{equation}
\PsiS{n}(0;\gamma^2) = -\pi/2 (n+1),
\label{introduction:xi2}
\end{equation}
and our expansions are (essentially) functions of the parameter
\begin{equation}
\xi = -\frac{2}{\pi} \PsiS{\nu}(0;\gamma^2)-1.
\label{introduction:xi}
\end{equation}
Among other things, this ensures that
$\chi(\xi;\gamma^2)$ is equal to the Sturm-Liouville eigenvalue
of $\ps{n}(x;\gamma^2)$ when $\xi = n$.  
In fact, for technical reasons we will discuss in detail when we describe
the numerical procedure used to construct them, these expansions are functions
of a variable related to $\xi$ through $\gamma$.  
However, for most intents
and purposes, they can be viewed as functions of $\xi$, and we will discuss
as if they are.
While $\chi$ is discontinuous as a function of characteristic
exponent, it is smooth as a function of $\xi$.  
This is rather dramatically demonstrated by  Figure~\ref{figure:chiplots}, 
which contains plots of $\chi$ as a function
of the characteristic exponent $\nu$ and of $\chi$ as a function
of the parameter $\xi$ when $\gamma=500$.
The other quantities we represent in this fashion
are the values at $0$ of the first three derivatives of the nonoscillatory phase
function $\PsiS{\nu}(x;\gamma^2)$ with respect to the argument $x$.  We will use the notations
\begin{equation}
\PsiS{\nu}'(x;\gamma^2), \ \ \PsiS{\nu}''(x;\gamma^2), \ \
\PsiS{\nu}'''(x;\gamma^2), \ \ \ldots
\end{equation}
to denote the derivatives of the nonoscillatory phase function  with respect to $x$.

Our  expansions of $\chi$ and the derivatives of the nonoscillatory phase
provide a mechanism for the  numerical evaluation of $\chi_n(\gamma^2)$ 
that runs  in time independent of the parameters $n$ and $\gamma$, and they
 give us  the data necessary to   calculate
 the nonoscillatory phase function $\PsiS{n}(x;\gamma^2)$
 representing $\ps{n}(x;\gamma^2)$.
The time required to construct $\PsiS{n}(x;\gamma^2)$ is independent
of $n$ and $\gamma$, as is the time required to evaluate
 $\ps{n}(x;\gamma^2)$ using $\PsiS{n}(x;\gamma^2)$.

The approach of \cite{BremerKummer} is designed for the 
regime in which the coefficients in (\ref{introduction:swe})
are sufficiently large and it loses accuracy when this is not the case.
Moreover, the use of precomputed expansions means that we must
{\it a priori} fix some range for the parameters
$\gamma$ and $n$.  The algorithm we describe here applies
when
\begin{equation}
 256 = 2^8 \leq \gamma \leq 2^{20} = 1,048,576  \ \ \mbox{and} \  \ 200 \leq n \leq \gamma.
\end{equation}
That our method doesn't apply when both $\gamma$ and $n$ are small is of little account
as the Xiao-Rokhlin algorithm is highly effective  in that
regime.  Moreover, our algorithm could be easily altered to allow
for the evaluation of $\ps{n}(x;\gamma^2)$ in the case of
larger values of $n$ and $\gamma$.
That the algorithm does not apply for small
$n$ and large $\gamma$, however, is a significant limitation.
The authors will report on an alternate method which can be used
in this  regime at a later date.

The remainder of this document is organized as follows.
In Section~\ref{section:preliminaries}, we carefully establish the notation we use
for the spheroidal wave function  and  review some well-known facts regarding them.
 Section~\ref{section:algorithm1} describe the method used to construct
the expansions of $\chi$ and of the quantities
\begin{equation}
\PsiS{\nu}'(0;\gamma^2), \ \ \PsiS{\nu}''(0;\gamma^2), \ \ \mbox{and} \ \
\PsiS{\nu}'''(0;\gamma^2)
\end{equation}
as functions of $\xi$ and $\gamma$.
  In Section~\ref{section:algorithm2}, we detail our algorithm for the
evaluation of $\ps{n}(x;\gamma^2)$.  Finally,  in Section~\ref{section:experiments},
we present the results of numerical experiments which demonstrate the
properties of our scheme.

\begin{section}{Preliminaries}
\label{section:preliminaries}

In this section, we set our notation for the spheroidal wave function
and briefly review certain well-known facts
which are used in the design of the algorithms of this paper.
We state without proof many assertions
regarding spheroidal wave functions.
We refer the reader to \cite{Meixner},  \cite{Imam}, \cite{Flammer},
and \cite{Arscott} for thorough and rigorous discussions of this material.



\begin{subsection}{Characteristic exponents}
\label{section:charexp}

The functions  (\ref{introduction:ps}) can be  distinguished from other solutions of 
(\ref{introduction:swe}) through their behavior at infinity as well as by the boundary conditions 
(\ref{introduction:bc}).   Indeed, $\ps{n}(x;\gamma^2)$ admits an expansion at infinity of the form
\begin{equation}
z^\nu \sum_{k=-\infty}^\infty d_k z^{2k}
\label{introduction:infexp}
\end{equation}
with $\nu= n$, and the prolate spheroidal wave functions of order $0$
are the only solutions of (\ref{introduction:swe}) of this type.
In general, for any values of the parameters $\gamma$ and $\chi$,  
(\ref{introduction:swe})
admits a solution which has an expansion of the form (\ref{introduction:infexp}) around infinity.
When $\nu$ is not a half-integer, there is a second independent solution which can
expanded as
\begin{equation}
z^{1-\nu} \sum_{k=-\infty}^\infty d_k z^{2k}
\label{introduction:infexp2}
\end{equation}
at infinity.  For half-integer values of $\nu$, the second independent solution
 takes  on the form
\begin{equation}
z^\nu\log(z)  \sum_{n=-\infty}^\infty d_n z^{2n} + z^\nu \sum_{n=-\infty}^\infty e_n z^{2n}
\end{equation}
at infinity.  The complex number $\nu$  in the expansion (\ref{introduction:infexp})
is called a characteristic exponent for the solutions of (\ref{introduction:swe}).
Obviously, $\nu$ is not uniquely determined by $\chi$ and $\gamma$.   

This ambiguity can be
resolved by considering what happens when $\gamma=0$.  In that event,
(\ref{introduction:swe}) becomes Legendre's differential,
and the relationship between $\chi$ and $\nu$ is well-known:
$\chi = \nu (\nu+1)$.  To define $\chi_\nu(\gamma^2)$ uniquely
for each $\gamma >0$ and $\nu$ which is not a half-integer, we require that
$\chi_\nu(\gamma^2)$ 
converge to $\nu(\nu+1)$ continuously as $\gamma \to 0^+$.  

The definition of $\chi_\nu(\gamma^2)$ in the case of half-integer
values of $\nu$ is more complicated. 
Indeed,  essentially all aspects of the spheroidal wave functions  
of half-integer orders require special attention.
Since we have no need of these functions, we will always assume
without further comment 
that $\nu$ is not a half-integer in what follows.
For a definitions of $\chi_\nu(\gamma^2)$ and the spheroidal
wave functions in the case of half-integer
values of $\nu$, we refer the reader to \cite{Meixner}, \cite{Meixner2}
and \cite{Imam}. 

\end{subsection}


\begin{subsection}{Notations for the spheroidal wave functions}

Several different notations for the spheroidal wave functions
are in widespread use.    For the most part, we follow Chapter~30 of \cite{DLMF}.    
In particular, we use 
$\PS{\nu}\left(z;\gamma^2\right)$
and
$\QS{\nu}\left(z;\gamma^2\right)$
to denote the angular spheroidal wave functions of the first and second kinds
of order $0$, characteristic exponent $\nu$ and bandlimit $\gamma$,
respectively.        The function $\PS{\nu}\left(z;\gamma^2\right)$
can defined via an expansion of the form
%
\begin{equation}
\sum_{k=-\infty}^\infty a_{\nu,k}(\gamma^2) P_{\nu+2k}(z),
\label{preliminaries:pexp}
\end{equation}
where $P_{\nu}$ is the Legendre function of the first kind of degree
$\nu$ and the coefficients $\{a_{\nu,n}(\gamma^2)\}$ are determined by a three-term
recurrence relation.
When $\nu$ is not an integer,
the function $\QS{\nu}\left(z;\gamma^2\right)$ admits the expansion
\begin{equation}
\sum_{k=-\infty}^\infty a_{\nu,k}(\gamma^2) Q_{\nu+2k}(z),
\label{preliminaries:qexp}
\end{equation}
where $Q_\mu$ denotes the Legendre function of second kind of degree $\mu$
and  the $\{a_{\nu,n}\}$ are as in (\ref{preliminaries:pexp}).
When viewed as a function of $\mu$, $Q_\mu(z)$
has simple poles at the negative integers.  Consequently, the representation
(\ref{preliminaries:qexp}) is not viable when $\nu$ is an integer.  
However, for nonnegative integers $n$, $\QS{2n}(z;\gamma^2)$ can be represented in the form
\begin{equation}
\sum_{k=0}^\infty b_{2n,k}(\gamma^2) P_{2k+1}(z)
+
\sum_{k=0}^\infty a_{2n,k}(\gamma^2) Q_{2k}(z)
\label{preliminaries:qeven}
\end{equation}
with $\{a_{2n,k}(\gamma^2)\}$
as before 
 and $\{b_{2n,k}(\gamma^2)\}$ a second set of coefficients which can be obtained
from $\{a_{2n,k}(\gamma^2)\}$ by solving a system of linear algebraic equations.
Likewise,  $\QS{2n+1}(z;\gamma^2)$ can be represented as a sum of the form
\begin{equation}
\sum_{k=0}^\infty b_{2n+1,k}(\gamma^2) P_{2k}(z)
+
\sum_{k=0}^\infty a_{2n+1,k}(\gamma^2) Q_{2k+1}(z).
\label{preliminaries:qodd}
\end{equation}

The function of the first kind $\PS{\nu}(z;\gamma^2)$ is analytic
on the cut plane  $\mathbb{C}\setminus\left(-\infty,-1\right]$,
and that of the second kind
 $\QS{\nu}(z;\gamma^2)$ is analytic on 
 $\mathbb{C}\setminus\left(-\infty,1\right]$.

The standard real-valued solutions of (\ref{introduction:swe})
on the cut $(-1,1)$
are  defined for $-1 <x <1$ via the formulas
\begin{equation}
\ps{\nu}\left(x;\gamma^2\right)
= \lim_{y\to 0} 
\PS{\nu}\left(x+iy;\gamma^2\right)
\end{equation}
and
\begin{equation}
\qs{\nu}(x;\gamma^2)
= \lim_{y\to 0^+} 
\frac{1}{2}
\left(
\QS{\nu}(x+iy;\gamma^2)
+
\QS{\nu}(x-iy;\gamma^2)
\right).
\end{equation}
They  are analogs of Ferrer's versions  $\p{\nu}(x)$ and $\q{\nu}(x)$ of the Legendre functions,
 which are the standard real-valued solutions of the Legendre's
differential equation defined on $(-1,1)$ (see, for instance, Chapters~14 of
\cite{DLMF}).     Clearly, $\ps{n}(x;\gamma^2)$ admits  the expansion
\begin{equation}
\ps{\nu}(x;\gamma^2) = \sum_{k=-\infty}^\infty a_{\nu,k}(\gamma^2) \p{\nu+2k}(z),
\label{preliminaries:pexp2}
\end{equation}
and  for $\nu$ which are not integers, we have
\begin{equation}
\qs{\nu}(x;\gamma^2) = \sum_{k=-\infty}^\infty a_{\nu,k}(\gamma^2) \q{\nu+2k}(z).
\label{preliminaries:qexp2}
\end{equation}
Moreover, the obvious analogs of the expansions (\ref{preliminaries:qeven})
and (\ref{preliminaries:qodd}) for $\qs{2n}(x;\gamma^2)$
and $\qs{2n+1}(x;\gamma^2)$ also hold.
The connection formula 
\begin{equation}
\lim_{y\to 0^+} \QS{\nu}(x+iy;\gamma^2) = 
\qs{\nu}(x;\gamma^2)
-i \frac{\pi}{2} \ps{\nu}(x;\gamma^2)
\label{preliminaries:qconnection}
\end{equation}
follows readily from the analogous formula for Legendre functions
(which can be found, for instance, in Chapter~14 of \cite{DLMF}).

We denote the radial spheroidal wave functions of the first
and second kinds of bandlimit $\gamma$, 
characteristic exponent $\nu$ and order $0$ via
$S^{(1)}_\nu(z;\gamma^2)$
and
$S^{(2)}_\nu(z;\gamma^2)$, respectively.  For values of $\nu$ which are
not half-integers, they admit expansions of the form
\begin{equation}
S^{(1)}_\nu(z;\gamma^2)=
\frac{1}{A_{\nu}(\gamma^2)} 
\sum_{k=-\infty}^\infty (-1)^k a_{\nu,k}(\gamma^2) \ \sqrt{ \frac{\pi}{2\gamma z} } J_{\nu+\frac{1}{2}+2k}(\gamma z),
\end{equation}
and
\begin{equation}
S^{(2)}_\nu(z;\gamma^2)=
\frac{1}{A_{\nu}(\gamma^2)} 
\sum_{k=-\infty}^\infty (-1)^k a_{\nu,k}(\gamma^2) \ \sqrt{ \frac{\pi}{2\gamma z} } Y_{\nu+\frac{1}{2}+2k}(\gamma z),
\end{equation}
where $J_\mu$ and $Y_\mu$ denote the Bessel function of the first and second kinds of order $\mu$,
respectively.   The coefficients $\{a_{\nu,k}(\gamma^2)\}$ are as in (\ref{preliminaries:pexp}) and
the normalizing constant $A_\nu(\gamma^2)$ is defined via
\begin{equation}
A_\nu(\gamma^2) = \sum_{k=-\infty}^\infty a_{\nu,k}.
\end{equation}
The radial spheroidal wave function of the third kind of bandlimit $\gamma$,
characteristic exponent $\nu$ and order $0$ is 
\begin{equation}
S^{(3)}_\nu(x;\gamma^2) = S^{(1)}_\nu(x;\gamma^2) + i S^{(3)}_\nu(x;\gamma^2),
\end{equation}
and it admits the expansion
\begin{equation}
S^{(3)}_\nu(z;\gamma^2)=
\frac{1}{A_{\nu}(\gamma^2)} 
\sum_{k=-\infty}^\infty (-1)^k a_{\nu,k}(\gamma^2) \ 
\sqrt{ \frac{\pi}{2\gamma z} } H^{(1)}_{\nu+\frac{1}{2}+2k}(\gamma z)
\label{preliminaries:S3exp}
\end{equation}
with $H^{(1)}_\mu(z)$ the Hankel function of the first kind of order $\mu$.
The radial spheroidal wave functions are analytic in the cut plane
$\mathbb{C}\setminus\left(-\infty,0\right]$.

It is in our notation for the Sturm-Liouville eigenvalues that we 
deviate from \cite{DLMF}.  There, 
 $\lambda_\nu(\gamma^2)$ is used to denote the Sturm-Liouville
eigenvalue of $\ps{\nu}(x;\gamma^2)$ with respect to the operator
\begin{equation}
\widetilde{L}_\gamma\left[y\right](x) = - (1-x^2)y'(x) + 2 x y'(x) - \gamma^2 (1-x^2) y(x),
\end{equation}
whereas we use  $\chi_{\nu}(\gamma^2)$ to denote the Sturm-Liouville eigenvalue
of $\ps{\nu}(x;\gamma^2)$ with respect to the operator $L_\gamma$ defined in
 (\ref{introduction:diffop}).
Obviously, $\chi_\nu(\gamma^2)$ is related to $\lambda_\nu(\gamma^2)$ via
\begin{equation}
 \chi_\nu(\gamma^2)  = \lambda_\nu(\gamma^2) + \gamma^2.
\end{equation}
Our convention is consistent with \cite{Osipov-Rokhlin-Xiao} and \cite{Flammer}.


\end{subsection}

\begin{subsection}{Phase functions for second order differential equations}

We say that a smooth function $\alpha$  is a phase function
for a second order differential equation of the form
\begin{equation}
y''(x) + q(x) y(x) = 0 \ \ \ \mbox{for all}\ \ a< x <b
\label{phase:ode}
\end{equation}
provided $\alpha'(x) > 0$ for $a < x <b$ and the functions
\begin{equation}
u(x) = \frac{\sin\left(\alpha(x)\right)}{\sqrt{\alpha'(x)}}
\label{phase:u}
\end{equation}
and
\begin{equation}
v(x) = \frac{\cos\left(\alpha(x)\right)}{\sqrt{\alpha'(x)}}
\label{phase:v}
\end{equation}
constitute a basis in the space of solutions of (\ref{phase:ode}).
Any second order differential equation can be converted into the form (\ref{phase:ode})
via a simple transformations.  For instance,
if $y$ satisfies (\ref{introduction:swe}) on $-1 < x < 1$, then
\begin{equation}
\varphi(x) = y(x) \sqrt{1-x^2}
\end{equation}
solves 
\begin{equation}
\varphi''(x) + q(x)\varphi(x) = 0
\ \ \ \mbox{for all}  \ \ -1 < x <1,
\label{phase:swe}
\end{equation}
with $q$ given by 
\begin{equation}
q(x) = \frac{1}{(1-x^2)^2}
+ \frac{\chi - \gamma^2 x^2}{1-x^2}.
\label{phase:q})
\end{equation}
We refer to (\ref{phase:swe}) as the normal form of the spheroidal wave equation.

Any pair of real-valued solutions of (\ref{phase:ode})
whose Wronskian is $1$
determines a phase function for (\ref{phase:ode}) up to 
a constant multiple of $2\pi$.
Indeed, (\ref{phase:u}) and (\ref{phase:v}) immediately
imply
\begin{equation}
\alpha'(t) = \frac{1}{(u(t))^2 + (v(t))^2},
\label{phase:ap}
\end{equation}
which determines $\alpha$ up to a constant and that 
constant is fixed modulo $2\pi$  by the requirement
that   (\ref{phase:u}) and (\ref{phase:v}) hold.

\end{subsection}


\begin{subsection}{Connection formulas for the radial spheroidal wave functions}

By examining the series expansions of the spheroidal wave functions
and  the angular wave functions at infinity, formulas
connecting the two can be obtained.    
Of particular interest to us are connection formulas for the 
radial spheroidal wave functions of the third kind of
integer characteristic exponents.  The boundary values
of these functions on the real line give rise to 
a pair of  real-valued solutions of   (\ref{introduction:swe})
which generate the nonoscillatory phase functions 
$\PsiS{n}(x;\gamma^2)$ we use
to represent the prolate spheroidal wave functions
$\ps{n}(x;\gamma^2)$.

In the case of nonnegative even integer characteristic exponents
we have
\begin{equation}
S^{(1)}_{2n}(z;\gamma^2) = K_{2n}^{(1)}(\gamma^2) \PS{2n}(z;\gamma^2)
\end{equation}
and
\begin{equation}
S^{(2)}_{2n}(z;\gamma^2) = K_{2n}^{(2)}(\gamma^2) \QS{2n}(z;\gamma^2),
\end{equation}
where
\begin{equation}
K_{2n}^{(1)}(\gamma^2) 
= \frac{(-1)^n }{A_{2n}(\gamma^2)} \frac{a_{2n,-n}(\gamma^2)}{\ps{2n}(0;\gamma^2)}
\end{equation}
and
\begin{equation}
K_{2n}^{(2)}(\gamma^2) 
=  \frac{(-1)^{n+1}}{\gamma A_{2n}(\gamma^2)} \frac{\ps{2n}(0;\gamma^2)}{a_{2n,-n}(\gamma^2)}.
\end{equation}
In the case of  nonnegative odd integer characteristic exponents,
\begin{equation}
S^{(1)}_{2n+1}(z;\gamma^2) = K_{2n+1}^{(1)}(\gamma^2) \PS{2n+1}(z;\gamma^2)
\end{equation}
and
\begin{equation}
S^{(2)}_{2n+1}(z;\gamma^2) = K_{2n+1}^{(2)}(\gamma^2) \QS{2n+1}(z;\gamma^2),
\end{equation}
where
\begin{equation}
K_{2n+1}^{(1)}(\gamma^2) 
= \frac{(-1)^n}{A_{2n+1}(\gamma^2)} 
\frac{\gamma}{3} \frac{a_{2n+1,-n}(\gamma^2)}{\ps{2n+1}'(0;\gamma^2)}
\end{equation}
and
\begin{equation}
K_{2n+1}^{(2)}(\gamma^2) 
=  \frac{(-1)^{n+1}}{\gamma A_{2n}(\gamma^2)} \frac{3}{\gamma^2} \frac{\ps{2n+1}'(0;\gamma^2)}
{a_{2n+1,-n}(\gamma^2)}.
\end{equation}
Here, we are using the convention that the prime
symbol indicates differentiation with respect to the
argument $x$ so that
$\ps{2n+1}'(0;\gamma^2)$ denotes the value of the derivative with respect to 
$x$ of $\ps{2n+1}(x;\gamma^2)$ at $0$.   For any nonnegative
integer value of $n$ we have, by virtue of the preceding formulas
and (\ref{preliminaries:qconnection}),
\begin{equation}
\lim_{y\to 0^+} S^{(3)}_n(x+iy;\gamma^2) =
D_{n}^{(1)}(\gamma^2) \ps{n}(x;\gamma^2)
+
i D_{n}^{(2)}(\gamma^2) \qs{n}(x;\gamma^2),
\end{equation}
where 
\begin{equation}
D_{n}^{(1)}(\gamma^2) = K_{n}^{(1)}(\gamma^2) + \frac{\pi}{2}K_{n}^{(2)}(\gamma^2) 
\label{phase:d1}
\end{equation}
and
\begin{equation}
D_{n}^{(2)}(\gamma^2) = K_{n}^{(2)}(\gamma^2).
\label{phase:d2}
\end{equation}
For noninteger values of $\nu$, there also exist coefficients  
$D_{\nu}^{(1)}(\gamma^2)$
and $D_{\nu}^{(2)}(\gamma^2)$ such that
\begin{equation}
\lim_{y\to 0^+} S^{(3)}_\nu(x+iy;\gamma^2) =
D_{\nu}^{(1)}(\gamma^2) \ps{\nu}(x;\gamma^2)
+
i D_{\nu}^{(2)}(\gamma^2) \qs{\nu}(x;\gamma^2).
\end{equation}
Their definitions, which are somewhat  more complicated than in the case
of integer characteristic exponents, 
can be found in Section~3.66 of  \cite{Meixner}.

\begin{figure}[b]
\begin{center}
\includegraphics[width=.49\textwidth]{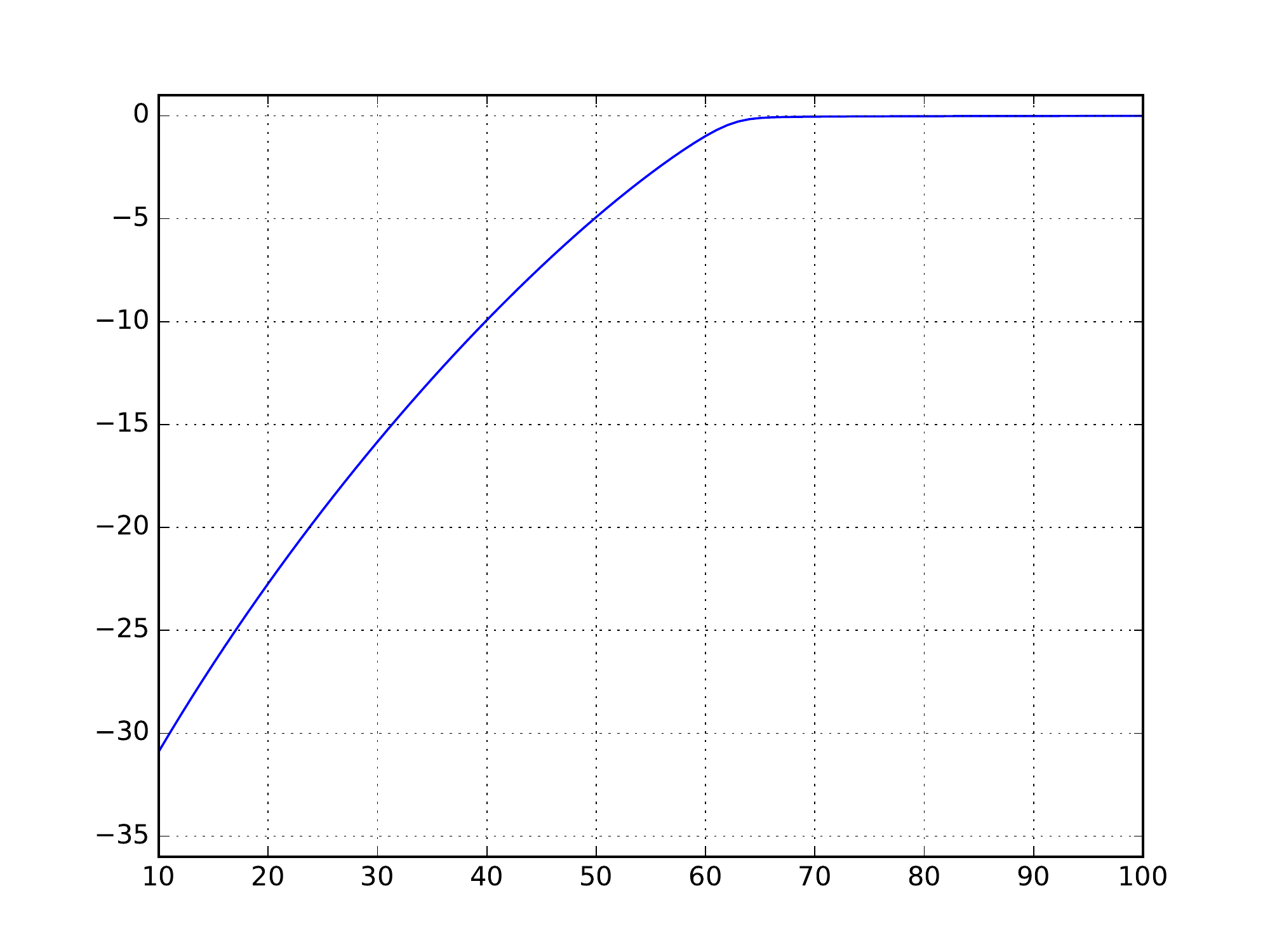}
\end{center}
\caption{A plot of the base-10 logarithm of $A_\nu(\gamma^2)$
as a function of $\nu$ when $\gamma=100$.}
\label{figure:wronskian}
\end{figure}

The Wronskian of the pair of solutions $\ps{\nu}(x;\gamma^2)$,
 $\qs{\nu}(x;\gamma^2)$ is 
\begin{equation}
\frac{A_\nu^2(\gamma^2)}{1-x^2}.
\end{equation}
When $n$ is small relative to $\gamma$, the magnitude of $A_\nu(\gamma^2)$ is extremely small.
See, for instance, 
Figure~\ref{figure:wronskian}, which  contains a plot of the base-10 logarithm of
$A_\nu(\gamma^2)$ as a function of $\nu$ when $\gamma=100$.  
Among other things it shows that when $\gamma=100$, the magnitude of $A_\nu(\gamma^2)$ 
already falls below $10^{-30}$.  The situation becomes even worse as $\gamma$ increases.
Clearly,   the pair $\ps{\nu}(x;\gamma^2)$ and $\qs{\nu}(x;\gamma^2)$ 
constitute a  basis  which is  extremely   ill-conditioned numerically for many
values of $\nu$ and $\gamma$.  This motivates the following definitions.  We let
\begin{equation}
C^{(1)}_\nu(\gamma^2) =  \frac{1}{A_\nu(\gamma^2)} D^{(1)}_\nu(\gamma^2)
\label{phase:c1}
\end{equation}
and
\begin{equation}
C^{(2)}_\nu(\gamma^2) =  \frac{1}{A_\nu(\gamma^2)} D^{(2)}_\nu(\gamma^2),
\label{phase:c2}
\end{equation}
so that 
\begin{equation}
u_\nu(x; \gamma^2) = C^{(1)}_\nu(\gamma^2) \ps{\nu}(x;\gamma^2)\sqrt{1-x^2}
\label{swe:u}
\end{equation}
and
\begin{equation}
v_\nu(x; \gamma^2) = C^{(2)}_\nu(\gamma^2) \qs{\nu}(x;\gamma^2)\sqrt{1-x^2}
\label{swe:v}
\end{equation}
is a pair of  solutions of  the normalized spheroidal wave equation (\ref{phase:swe})
whose  Wronskian is $1$.

\vskip 1em
\begin{remark}
The numerical evaluation of the coefficients 
$C^{(1)}_\nu(\gamma^2)$ and $C^{(2)}_\nu(\gamma^2)$
through the
formulas (\ref{phase:c1}), (\ref{phase:c2}) (\ref{phase:d1}) and
(\ref{phase:d2}) is problematic.
When $\gamma$ is of large magnitude $n$  is small relative to $\gamma$, 
the evaluation of these 
formulas using finite precision arithmetic results in 
 catastrophic  cancellation errors.  In fact, when $\gamma$ is of large
magnitude, this is the case even for relatively large values of $n$ 
(for instance, $\gamma = 1000$ and $n = 500$).  
We do not make use of these formulas in the algorithm of this paper.
\end{remark}

\end{subsection}


\begin{subsection}{The nonoscillatory phase function for the spheroidal wave equation}
\label{section:nonoscillatory}

It follows from the formula
\begin{equation}
S^{(3)}_\nu(z) = 
-\frac{\exp\left(-i\frac{\pi}{2}\nu\right)}{A_\nu(\gamma^2)}
\int_{1}^\infty \exp(i\gamma z t) \ps{\nu}(t;\gamma^2)\ dt,
\label{preliminaries:S3Fourier}
\end{equation}
which specifies the Fourier transform of the radial spheroidal wave function
of the third kind and can be found in Section~3.84 of \cite{Meixner},
that   the function
\begin{equation}
f(x) = \lim_{y\to0^+} \left|S^{(3)}_\nu(x+iy;\gamma^2)\right|^2 
\end{equation}
is absolutely monotone on the interval $(-1,1)$.  That is, 
$f(x)$ and its  derivatives of all orders are positive on $(-1,1)$.
Indeed, this result can be obtained by using (\ref{preliminaries:S3Fourier})
to derive a formula expression the  Laplace transform of 
the boundary value of 
\begin{equation}
 \left|S^{(3)}_\nu(x+iy;\gamma^2)\right|^2
\end{equation}
as a convolution of angular  spheroidal wave functions of the first kind.

\begin{figure}[t!!!]
\begin{center}
\includegraphics[width=.49\textwidth]{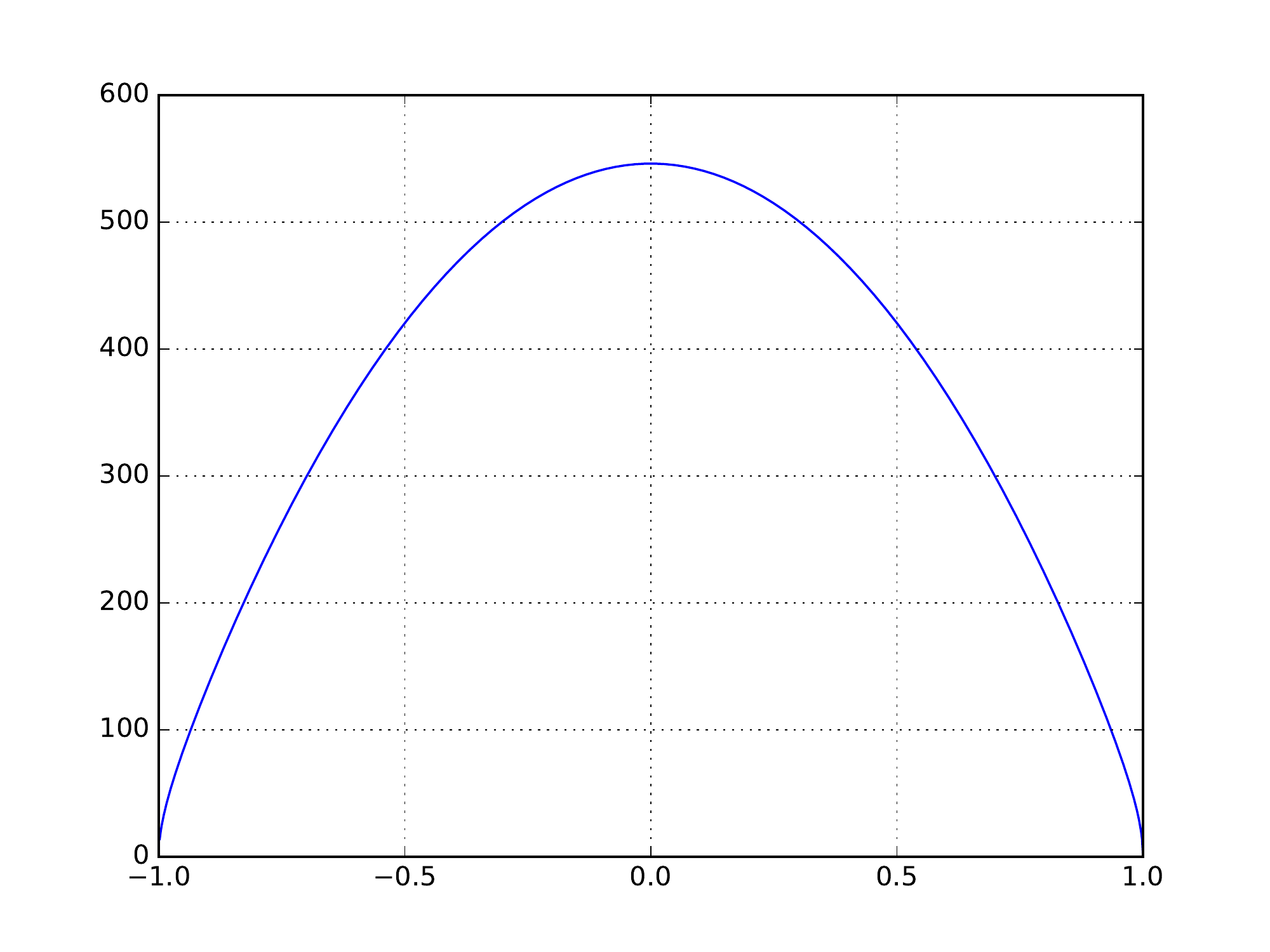}
\hfil
\includegraphics[width=.49\textwidth]{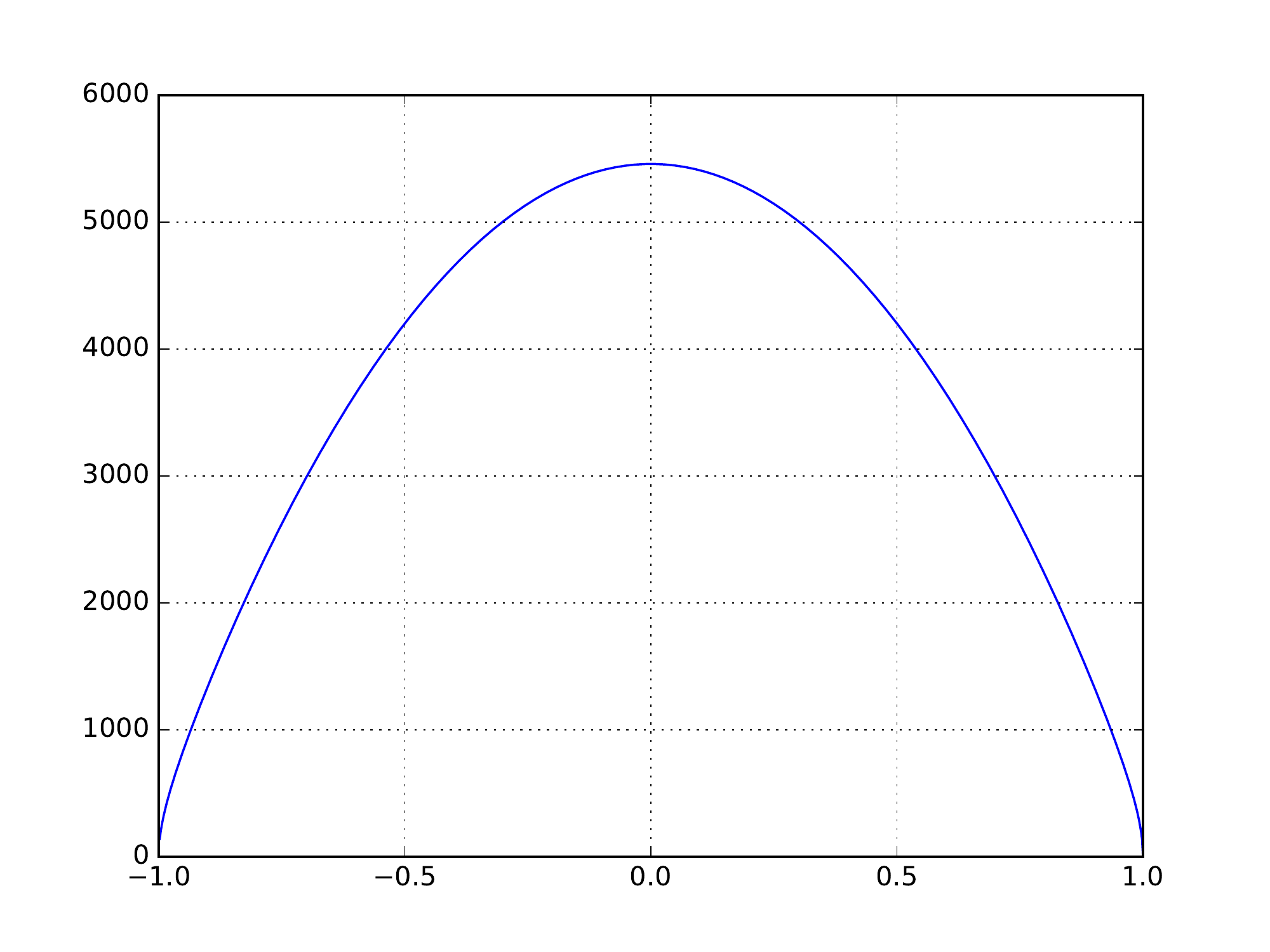}
\end{center}
\caption{On the left is a plot of $\PsiS{\nu}'(x;\gamma^2)$
when  $\gamma = 500$ and $\nu = 400$, and on the right
is a plot of  $\PsiS{\nu}'(x;\gamma^2)$ when $\gamma=5000$
and $\nu = 4000$.}
\label{figure:phaseder}
\end{figure}

We use $\PsiS{\nu}(x;\gamma^2)$ to denote a phase
function for the normal form of the spheroidal
wave equation (\ref{phase:swe}) which gives 
rise to the solutions
(\ref{swe:u}) and (\ref{swe:v}) 
via the formulas
\begin{equation}
u_{\nu}(x;\gamma^2) = \frac{\sin\left(\PsiS{\nu}(x;\gamma^2)\right)}{\sqrt{\PsiS{\nu}'(x;\gamma^2)}}
\end{equation}
and
\begin{equation}
v_{\nu}(x;\gamma^2) = \frac{\cos\left(\PsiS{\nu}(x;\gamma^2)\right)}{\sqrt{\PsiS{\nu}'(x;\gamma^2)}}.
\end{equation}
We uniquely determine $\PsiS{\nu}(x;\gamma^2)$ by requiring that
\begin{equation}
\lim_{x\to 1^-} \PsiS{\nu}(x;\gamma^2) = 0.
\end{equation}
We note that the derivative of $\PsiS{\nu}(x;\gamma^2)$ with respect
to $x$ is positive, so $\PsiS{\nu}(x;\gamma^2)$ is a negative function
which increases to $0$ as  $x\to 1$ from the left.
According to (\ref{phase:ap}),
\begin{equation}
\PsiS{\nu}'(x;\gamma^2) = 
\frac{1}
{
\left(u_\nu(x;\gamma^2)\right)^2
+
\left(v_\nu(x;\gamma^2)\right)^2
},
\label{phase:swe_ap}
\end{equation}
where we are once again using the convention that the prime symbol denotes differentiation
with respect to the argument $x$.  From (\ref{swe:u}) and (\ref{swe:v})
is is clear that the reciprocal of (\ref{phase:swe_ap})
\begin{equation}
W_\nu(x;\gamma^2) = \left(u_\nu(x;\gamma^2)\right)^2 + \left(v_\nu(x;\gamma^2)\right)^2
\label{phase:W}
\end{equation}
is a constant multiple of 
\begin{equation}
\left|\lim_{y\to0^+}S^{(3)}_\nu(x+iy;\gamma^2)\right|^2 (1-x^2).
\end{equation}
In particular, $\PsiS{\nu}(x;\gamma^2)$ is nonoscillatory in the sense
that the reciprocal of its derivative is equal to $(1-x^2)$
times an absolutely monotone function.    
This is an extremely strong notion of ``nonoscillatory,'' and while many second
order differential equations admit a phase function which 
is nonoscillatory in some sense, it is rare 
that they admit a  phase function which is related by a sequence of algebraic
operations to an absolutely monotone or completely monotone function.
See \cite{Bremer-Rokhlin} for a much more general notion
of nonoscillatory phase function which applies to a large
class of second order differential equations.
Figure~\ref{figure:phaseder} contains the plots of
the derivative of $\PsiS{\nu}(x;\gamma^2)$ for
two different pairs of the parameters $\gamma$ and $\nu$.

\vskip 1em
\begin{remark}
Formula (\ref{preliminaries:S3Fourier}) follows from and is an  analog of 
\begin{equation}
\sqrt{\frac{\pi}{2z}} H^{(1)}_{\nu+\frac{1}{2}}(z) = - \exp\left(-i\frac{\pi}{2} \nu\right) 
\int_1^\infty \exp(izt) \mathsf{P}_\nu(t)\ dt,
\label{preliminaries:hankelFourier}
\end{equation}
which specifies the Fourier transform of the spherical Hankel function
of the first kind of degree $\nu$.
From (\ref{preliminaries:hankelFourier}) and standard results
regarding the Fourier transform, it follows that
\begin{equation}
\left|\sqrt{\frac{\pi}{2z}} H^{(1)}_{\nu+\frac{1}{2}}(z)\right|^2 
= 
\int_0^\infty \exp(izt) \mathsf{P}_\nu\left(1+\frac{t^2}{2}\right)\ dt,
\end{equation}
which can be rearranged as
\begin{equation}
\frac{1}{z}J_{\nu+\frac{1}{2}}^2(z) + \frac{1}{z}Y_{\nu+\frac{1}{2}}^2(z)
= \frac{2}{\pi}
\int_0^\infty \exp(-zt) \mathsf{P}_\nu\left(1+\frac{t^2}{2}\right)\ dt.
\label{preliminaries:nicholson}
\end{equation}
%
Since $\mathsf{P}_\nu(t)$ is nonnegative on $(1,\infty)$, we have that
\begin{equation}
\frac{1}{z}J_{\nu+\frac{1}{2}}^2(z) + \frac{1}{z} Y_{\nu+\frac{1}{2}}^2(z)
\label{preliminaries:bessel}
\end{equation}
is completely monotone on the interval $(0,\infty)$.  A smooth function $f$ is completely
monotone on an interval $(a,b)$ if
\begin{equation}
(-1)^k f^{(k)}(x) \geq 0
\end{equation}
for $x \in (a,b)$ and all nonnegative integers $k$.  A function is completely
monotone on $(0,\infty)$ if and only if it is the Laplace transform
of a positive Borel measure.
The function (\ref{preliminaries:bessel}) is the reciprocal of the derivative
of a phase function for the normal form
\begin{equation}
y''(t) + \left(1  + \frac{\frac{1}{4}-\nu^2}{t^2}\right)y(t) = 0
\end{equation}
of Bessel's differential equation.  So the normal form of Bessel's
differential equation admits a phase function whose derivative
is the reciprocal of a completely monotone function.
We note that Formula~(\ref{preliminaries:nicholson}) is an analog of
Nicholson's classical integral representation formula
(see Section~13.73 of \cite{Watson}), which also implies
that (\ref{preliminaries:bessel}) is completely monotone.
\end{remark}









\end{subsection}


\begin{subsection}{Kummer's equation, Riccati's equation  and Appell's equation}

We now briefly discuss three differential equations which can be 
solved to calculate phase functions for (\ref{phase:ode}).
The second order nonlinear ordinary differential 
\begin{equation}
(\alpha'(x))^2 = q(x) - \frac{1}{2} \frac{\alpha'''(x)}{\alpha'(x)} 
+ \frac{3}{4} \left(\frac{\alpha''(x)}{\alpha'(x)}\right)^2
\label{eqs:kummer}
\end{equation}
satisfied by the derivative of phase
functions for (\ref{phase:ode})
can be obtained from (\ref{phase:ap}) through repeated differentiation.
We refer to (\ref{eqs:kummer}) as  Kummer's equation after E. E. Kummer who studied it in \cite{Kummer}.    
Kummer's equation can also be obtained by decomposing  the  Riccati equation
\begin{equation}
r'(x) + (r(x))^2 + q(x) = 0
\label{eqs:riccati}
\end{equation}
satisfied by the logarithmic derivatives of solutions of (\ref{phase:ode})
into real and imaginary parts.       
It can be verified through direct computation that if $u$ and $v$ are solutions of 
(\ref{phase:ode}) then  
\begin{equation}
w(x) = (u(x))^2 + (v(x))^2
\end{equation}
solves
\begin{equation}
w'''(x) + 4 q(x) w'(x) + 2 q'(x) w(t) = 0.
\label{eqs:appell}
\end{equation}
We refer to (\ref{eqs:appell}) as Appell's equation, after P. Appell
who discussed it in \cite{Appell}.

The scaled  phase function 
\begin{equation}
\PsiS{\nu}(x;\gamma^2)\sqrt{1-x^2}
\end{equation}
satisfies Kummer's equation (\ref{eqs:kummer}) with $q$ as in (\ref{phase:q}),
while the function $W_\nu(x;\gamma^2)$
defined via (\ref{phase:W}) satisfies Appell's equation (\ref{eqs:appell})
with $q$ as in (\ref{phase:q}).

\end{subsection}


\begin{subsection}{Chebyshev expansions}
An nth order univariate Chebyshev expansion on the interval $(a,b)$
is a sum of the form
\begin{equation}
\sum_{i=0 }^n \beta_{i}\ T_i\left(\frac{2}{b-a} x - \frac{b+a}{b-a}\right),
\label{chebyshev:uniexp}
\end{equation}
where $T_m(x) = \cos(m \arccos(x))$ is the Chebyshev polynomial of degree $m$.
We refer to the collection of points $t_0,t_1,\ldots,t_n$ defined by 
\begin{equation}
t_j =  \cos\left(\frac{j\pi}{n}\right), \ \ j=0,1,\ldots,n
\label{chebyshev:grid}
\end{equation}
as the nth order Chebyshev grid on the interval $[-1,1]$,
and the set of points
\begin{equation}
\frac{b-a}{2} t_j + \frac{b+a}{2}, \ \ j=0,1,\ldots,n
\label{chebyshev:grid2}
\end{equation}
as the (n+1)-point Chebyshev  grid on the interval $[a,b]$.
For any continuous function $f:[a,b] \to \mathbb{R}$,
we call  the unique expansion of the form (\ref{chebyshev:uniexp})
which agrees with $f$ at the nodes (\ref{chebyshev:grid2}) 
the nth order Chebyshev expansion of $f$ on $[a,b]$.
When $f$ is infinitely differentiable,
the nth order Chebyshev expansion of $f$ on $[a,b]$ converges to $f$ in the 
$C(\left[a,b\right])$ norm  superalgebraically  as n increases, and 
it converges to $f$ exponentially fast if 
$f$ is analytic in neighborhood of the interval $[a,b]$.
The widespread use of Chebyshev expansions
(and expansions in other families of orthogonal polynomials) in numerical calculations 
is principally due to their  favorable stability  properties.
The coefficients in the  Chebyshev expansion of $f$ on $[a,b]$ can be computed
in a numerically stable fashion from the values of $f$
at the nodes of the Chebyshev grid on $[a,b]$, and
 (\ref{chebyshev:uniexp}) is well-conditioned as a function of 
the coefficients $\{\beta_i\}$.  
We refer the  reader to \cite{Trefethen} for a thorough treatment of 
these and other related results in approximation theory.

An nth order bivariate Chebyshev expansion on the rectangle
$(a,b) \times (c,d)$  is a sum of the form
\begin{equation}
\sum_{0 \leq i+j \leq n} \beta_{i,j} T_i\left(\frac{2}{b-a} x - \frac{b+a}{b-a}\right) 
T_j\left(\frac{2}{d-c} y - \frac{d+c}{d-c}\right),
\label{chebyshev:biexp}
\end{equation}
and we call the collection of points
\begin{equation}
\left(
\frac{b-a}{2} t_i + \frac{b+a}{2},
\frac{d-c}{2} t_j + \frac{d+c}{2} 
\right),\ \ \ i,j=0,1,\ldots,m
\label{chebyshev:tensor_grid}
\end{equation}
where $t_0,t_1,\ldots,t_n$ are as in (\ref{chebyshev:grid}),
 the nth order Chebyshev  grid on the rectangle $(a,b) \times (c,d)$. 
For any continuous function $f:[a,b] \to \mathbb{R}$,
we call the unique expansion of the form (\ref{chebyshev:biexp})
which agrees with $f$ at the nodes (\ref{chebyshev:tensor_grid}) the
nth order bivariate Chebyshev expansion of $f$ on $[a,b]$.
The coefficients in such an expansion can be computed in a numerical
stable fashion from the values of $f$ at the nodes (\ref{chebyshev:tensor_grid}),
and the expansion (\ref{chebyshev:biexp}) is well-conditioned
as a function of its coefficients.
As in the case of univariate Chebyshev expansions, the $n$th order
bivariate Chebyshev expansion of an infinitely differentiable function $f$ converges
to $f$ superalgrebraically with increasing $n$, and 
the  analyticity of $f$ in a neighborhood of implies exponential convergence.

The  nth order piecewise Chebyshev expansion of the continuous
function $f:[a,b] \to \mathbb{R}$  
with respect to the partition
\begin{equation}
a= a_1 < a_2 < \ldots < a_m  = b
\label{chebyshev:intervals}
\end{equation}
of $[a,b]$ consists  of the nth order Chebyshev expansions of $f$  on each
of the intervals
\begin{equation}
(a_1,a_2), (a_2,a_3), \ldots, (a_{m-1},a_{m}).
\end{equation}
The nth order piecewise bivariate Chebyshev expansion of the continuous
function $f:[a,b] \times [c,d] \to \mathbb{R}$   with respect to the partitions 
\begin{equation}
a= a_1 < a_2 < \ldots < a_{m_1}  = b
\label{chebyshev:intervals2}
\end{equation}
and
\begin{equation}
c= c_1 < c_2 < \ldots < c_{m_1}  = d
\label{chebyshev:intervals3}
\end{equation}
consists  of the nth order bivariate  Chebyshev expansions of $f$  on each
of the rectangles
\begin{equation}
[a_i,a_{i+1}] \times [c_j,c_{j+1}], \ \ i=0,1,\ldots,m_1-1,\ \ j=0,1,\ldots,m_2-1.
\end{equation}
We generally prefer the use of piecewise expansions to a single high order
expansion for two reasons: they are more flexible in that a larger
class of functions (including many singular functions) can be represented
efficiently  using piecewise expansions, and, perhaps more importantly for this work,
the cost of evaluating a piecewise expansion at a single point is generally 
much lower.

\end{subsection}


\begin{subsection}{Adaptive Chebyshev Discretization}
\label{section:adap}

We now briefly describe a fairly standard procedure for
 adaptively discretizing a smooth function $f:[a,b] \to \mathbb{R}$.
It takes as input a desired precision $\epsilon >0$, a positive integer $n$
and a subroutine for evaluating $f$.
The goal of this procedure is to construct a partition
\begin{equation}
a = a_1 < a_2 < \cdots < \gamma_{m} = b
\end{equation}
of $[a,b]$ such that  the $n$th order Chebyshev expansion of $f$ on each of the subintervals
$[a_j,a_{j+1}]$ approximates $f$ with accuracy $\epsilon$.  That is,
for each $j=1,\ldots,m-1$ we aim to achieve
\begin{equation}
\sup_{x \in [a_j,a_{j+1}]}
\left| 
f(x) - 
\sum_{i=0}^n \beta_{i,j} T_{i} \left( \frac{2}{a_{j+1}-a_j } x + 
\frac{a_{j+1}+a_j}{a_j-a{j+1}} \right)
\right| < \epsilon,
\label{preliminaries:adaptive:1}
\end{equation}
where $\beta_{0,j},b_{1,j}\ldots,\beta_{n,j}$ are the coefficients in the $n$th order Chebyshev expansion
of $f$ on the interval $\left[a_j,a_{j+1}\right]$.  

During the procedure,  two lists of subintervals are maintained: a list
of subintervals which are to be processed and a list of output subintervals.
Initially, the list of subintervals to be processed consists of $[a,b]$ 
and the list of output subintervals is empty.
The procedure terminates when the list of subintervals to be processed is empty
or when the number of subintervals in this list  exceeds a present
limit (we usually take this limit to be $300$).
In the latter case, the procedure is deemed to have failed.
As long as the list of subintervals to process is nonempty and its length
does not exceed the preset maximum, the algorithm
proceeds by removing a subinterval $\left[\eta_1,\eta_2\right]$
from that list  and performing the following operations:
\begin{enumerate}

\item
Compute the coefficients $\beta_0,\ldots,\beta_n$ 
in the nth order Chebyshev expansion
of the restriction of $f$ on the $[\eta_1,\eta_2]$.

\vskip 1em
\item
Compute the quantity 
\begin{equation}
\Delta = 
\frac{\max\left\{
\left|\beta_{\frac{n}{2}+1} \right|, 
\left|\beta_{\frac{n}{2}+2} \right|, \ldots,
\left|\beta_n \right|
\right\}}
{\max\left\{
\left|\beta_0 \right|, 
\left|\beta_1 \right|, \ldots,
\left|\beta_n \right|
\right\}}.
\end{equation}

\vskip 1em
\item
If $\Delta < \epsilon$ then 
 the subinterval $\left[\eta_1,\eta_2\right]$ is added to the list
of output subintervals.

\vskip 1em
\item 

If $\Delta \geq \epsilon$, then the subintervals
\begin{equation}
\left[\eta_1,\frac{\eta_1+\eta_2}{2}\right]
\ \ \mbox{and}\ \ 
\left[\frac{\eta_1+\eta_2}{2}, \eta_2 \right]
\end{equation}
are added to the list of  subintervals to be processed.

\end{enumerate}

This algorithm is  heuristic in the sense that there is no guarantee
that (\ref{preliminaries:adaptive:1}) will be achieved, but 
similar adaptive discretization procedures are widely used
with great success. 

There is one common circumstance which 
leads to the failure of this procedure.  The quantity $\Delta$ is an attempt to estimate
the relative accuracy with which the Chebyshev expansion
of $f$ on the interval $\left[\eta_1,\eta_2\right]$ approximates
$f$.  In cases in which the condition number of the evaluation of
$f$  --- whose value at the point $x$ is 
\begin{equation}
\left|  x  \frac{f'(x)}{f(x)} \right|
\end{equation}
--- is larger than $\epsilon$ on some part of $[a,b]$,
the procedure will generally fail or an excessive number of subintervals
will be generated.   
Particular care needs to be taken when $f$ has a zero
in $[a,b]$.
In most cases, for $x$ near  a zero of $f$, the condition number of evaluation
of $f(x)$  
is large.   In this article, we avoid such difficulties  by only
applying this procedure to 
 functions which are bounded away from $0$.

\end{subsection}


\end{section}

\begin{section}{Numerical construction of the expansions of $\chi$ and the values
of the derivatives of the nonoscillatory phase function at $0$}
\label{section:algorithm1}

In this section, we describe the method which was used to construct the expansions of $\chi$ 
and the values
\begin{equation}
\PsiS{\nu}'(0;\gamma^2), \ \ \PsiS{\nu}''(0;\gamma^2), \ \ \mbox{and} \ \
\PsiS{\nu}'''(0;\gamma^2)
\label{alg1:functions}
\end{equation}
of the first few derivatives of the nonoscillatory phase function
at $0$.  
Our expansions are functions of $\gamma$ and a parameter
which is closely related to the quantity
$\xi$ defined via (\ref{introduction:xi}).
They take the form of   bivariate Chebyshev expansions
of order $k=29$.   After their construction,
they were written to a Fortran file on the disk for 
later use by the algorithm of Section~\ref{section:algorithm2};
each of them expansions occupies approximately $1.1$ megabyte of memory.
The computations described here were carried out on workstation
equipped with  $28$ Intel Xeon E5-2697 processor cores
running at 2.6 GHz.  They took approximately $20$ minutes to complete.

Our procedure made extensive use 
of the algorithm of \cite{BremerKummer},
which  allowed us to calculate the nonoscillatory phase function $\PsiS{\nu}(x;\gamma^2)$
and its first few derivatives given $\gamma$ and the value of $\chi_\nu(\gamma^2)$.
In particular, we used it as a mechanism for evaluating
\begin{equation}
\PsiS{\nu}(0;\gamma^2),\ \ \PsiS{\nu}'(0;\gamma^2), \ \ \PsiS{\nu}''(0;\gamma^2)
\ \ \mbox{and} \ \ \PsiS{\nu}'''(0;\gamma^2)
\label{alg1:functions2}
\end{equation}
as functions of $\gamma$ and $\chi$.

Our procedure began by introducing the partition
\begin{equation}
2^8 < 2^9 < 2^{10} < \ldots < 2^{18} < 2^{19} < 2^{20}
\label{alg1:gamma_part}
\end{equation}
of the interval 
\begin{equation}
256 = 2^8 \leq \gamma \leq 2^{20} = 1,048,576.
\label{alg1:gamma_int}
\end{equation}
We then formed the $(k+1)$-point Chebyshev grid on each of the intervals
defined by this partition.  For each  $\gamma$ in the resulting collection of points, 
we performed the following sequence of operations:

\begin{enumerate}
\item
We adaptively discretized the functions (\ref{alg1:functions2})  with respect
to the variable $\chi$ (with $\gamma$ held constant) over the interval
\begin{equation}
 \chi_{\xi_1}(\gamma^2) \leq \chi \leq \chi_{\xi_2}(\gamma^2),
\label{alg1:interval2}
\end{equation}
where $\xi_1 = 200$ and $\xi_2 = \gamma$.
Recall, that our expansions are meant to apply in the the
case of values of the parameter $\xi$ defined via (\ref{introduction:xi})
between $200$ and $\gamma$.
The scheme of Section~\ref{section:adap}
was used to perform this task; the order for the Chebyshev expansions
was taken to be $k$ and the requested precision was $\epsilon = 10^{-14}$.
The result  was a partition 
\begin{equation}
\eta_1 < \eta_2  < \ldots < \eta_m
\label{alg1:partition2}
\end{equation}
of  (\ref{alg1:interval2}), 
and the $k$th order piecewise Chebyshev expansions
of the  functions listed in   (\ref{alg1:functions2})
with respect to this partition.  We refer 
to the expansion of the value of the phase function via $\alpha(\chi)$,
the expansion of its second derivative via $\alpha'(\chi)$,
and so on.

\item
We next defined a function $\xi(\chi)$ via
\begin{equation}
\xi(\chi) = -\frac{2}{\pi}\alpha(\chi) -1.
\end{equation}
Because of (\ref{introduction:xi2}), the image of the interval
(\ref{alg1:interval2}) under this mapping is $[\xi_1,\xi_2]$.
We next formed the partition
\begin{equation}
\sigma_1 < \sigma_2 < \ldots < \sigma_m
\label{alg1:partition3}
\end{equation}
of $[\xi_1,\xi_2]$
%
%
by letting
\begin{equation}
\sigma_i = \xi(\eta_i),
\end{equation}
and  constructed the $k$th order piecewise Chebyshev expansion
of the inverse function $\chi(\xi)$ of $\xi(\chi)$
with respect to this partition.
We did so by computing the value of $\chi$ at each Chebyshev node
via the most  primitive root-finding method imaginable:  bisection.
The value of  $\chi$ increases monotonically with
increasing $\xi$, which made these computations significantly
simpler.

The inverse function of a polynomial of degree $k$ obviously
need not be a polynomial of degree $k$, and so
the   piecewise Chebyshev expansion of the inverse function
produced by a procedure of this sort can fail to accurately represent
it, even if the piecewise Chebyshev expansion of the original
function is highly accurate.  We relied on the facts
that the functions being inverted are extremely smooth,
 and that the discretizations formed
by the procedure of Section~\ref{section:adap} are somewhat
oversampled.  Moreover, we carefully verified the expansions of the inverse
functions generated in this step.

\item
%
We then defined a new parameter $\zeta$ (whose role will be made clear
shortly via
\begin{equation}
\xi = \xi_1 + (\xi_2 -\xi_1) \zeta,
\label{alg1:zeta}
\end{equation}
so that as $\zeta$ ranges over $(0,1)$, $\xi$ ranges 
over $(\xi_1,\xi_2)$.  We  introduced the partition
\begin{equation}
\zeta_1 < \zeta_2 < \ldots < \zeta_m
\label{alg1:partition4}
\end{equation}
of the interval $(0,1)$ which corresponds to (\ref{alg1:partition3})
and formed $k$th order piecewise Chebyshev expansions of the
functions
\begin{equation}
\chi(\zeta),\ \  \alpha'(\zeta),\ \  \alpha''(\zeta), \ \mbox{and} \ \ \alpha''(\zeta)
\end{equation}
with respect to (\ref{alg1:partition4}).
This can be done easily using the expansion of $\chi(\xi)$
formed in the preceding step of this procedure and the expansions
of these functions with respect to $\chi$ formed in the 
first step of this procedure.
\end{enumerate}
At this stage, for each point $\gamma$ which is a node in
one of the  $(k+1)$-point Chebyshev grids
on the intervals (\ref{alg1:gamma_part}),
we had 
 a partition (\ref{alg1:partition4})
and
 $k$th order Chebyshev expansions of $\chi$ and the quantities 
(\ref{alg1:functions2}) with respect to 
this partition.   The Chebyshev
expansion were functions of $\zeta$ and the partition
is of the interval $(0,1)$ over which $\zeta$ varies.

Next, we formed a single unified partition 
\begin{equation}
\tilde{\zeta}_1 < \tilde{\zeta}_2 < \ldots < \tilde{\zeta}_l
\label{alg1:joint}
\end{equation}
of $(0,1)$  by applying the adaptive procedure of Section~\ref{section:adap}
repeatedly to each of these expansions.  That is, we applied it to
the first expansion, and then used the resulting collection of intervals
as input while applying the procedure to the second expansion, and so on.
The requested precision for the discretization procedure was $\epsilon = 10^{-14}$.
The result was a collection of intervals sufficiently dense to discretize
each of the expansions, independent of $\gamma$.

We now had the ability to evaluate $\chi$ and the values at $0$
of the first three derivatives of the nonoscillatory phase functions
as functions of the parameter $\zeta$ for each value of $\gamma$
in  one of the  $(k+1)$-point Chebyshev grids
on the intervals (\ref{alg1:gamma_part}).  
This allowed us to form the $k$th order bivariate Chebyshev expansions
of these quantities with respect to the partitions
(\ref{alg1:gamma_part}) and (\ref{alg1:joint}).
These were the final product of the procedure of this section,
and the expansions which we use in the algorithm
of the following section.
We note that the value of $\xi_2$  in the relation (\ref{alg1:zeta}) 
defining $\zeta$ depends on $\gamma$.  


\end{section}

\begin{section}{An algorithm for the numerical Calculation of $\ps{n}(x;\gamma^2)$}
\label{section:algorithm2}

In this section, we describe our algorithm for the numerical
evaluation of $\ps{n}(x;\gamma^2)$.   It is divided into
two stages: a precomputation stage in which a piecewise
Chebyshev expansion of the nonoscillatory
phase function $\PsiS{n}(x;\gamma^2)$ is constructed, 
and an evaluation phase in which the phase function
is used to evaluate $\ps{n}(x;\gamma^2)$ at 
one or more points.
Owing to the symmetry of the functions $\ps{n}(x;\gamma^2)$ (they are even
functions when $n$ is even and odd functions when $n$ is odd), 
it is only necessary
to construct an expansion of 
$\PsiS{n}(x;\gamma^2)$ over the interval $[0,1)$.

The precomputation phase of the algorithm takes as input $n$ and $\gamma$.
We let
\begin{equation}
\zeta = \frac{n - \xi_1}{\xi_2 -\xi_1},
\end{equation}
where
\begin{equation}
\xi_1 = 200 \ \ \mbox{and} \ \ \xi_2 = \gamma.
\end{equation}
We next evaluate the precomputed expansions discussed in Section~\ref{section:algorithm1},
which are functions of $\zeta$ and $\gamma$,
to obtain the values of 
\begin{equation}
\chi_n(\gamma^2), \ \ \PsiS{n}(0;\gamma^2), \ \ \PsiS{n}'(0;\gamma^2), \ \
\PsiS{n}''(0;\gamma^2), \ \ \mbox{and}\ \ \PsiS{n}'''(0;\gamma^2).
\label{alg2:initial}
\end{equation}
The cost of evaluating these expansions is independent of $\gamma$ and $n$.

At this stage, we could solve an initial value problem for 
the differential equation (\ref{eqs:kummer})
to compute $\PsiS{n}(x;\gamma^2)$ on the interval $(0,1)$ --- this
is similar to the approach in 
 \cite{BremerKummer}, which operates by solving Kummer's equation.
However, for most values of $n$ and $\gamma$, the
spheroidal wave equation has turning points in 
the interval $(0,1)$, and the numerical solution of Kummer's
equation is complicated by the presence of turning points.
Instead of solving Kummer's equation to construct
 $\PsiS{n}(x;\gamma^2)$, we solve Appell's equation
 (\ref{eqs:appell})
to obtain the function   $W_n(x;\gamma^2)$ defined via (\ref{phase:W}).
As discussed in Section~\ref{section:nonoscillatory},
the function 
\begin{equation}
\frac{W_n(x;\gamma^2) }{(1-x^2)}
\label{alg2:nonosc}
\end{equation}
is absolutely monotone on the interval $(-1,1)$,
and the numerical solution of Appell's equation is not made
more difficult by the presence of turning points.
We use the quantities in (\ref{alg2:initial}) to compute
the values of 
\begin{equation}
W_n(0;\gamma^2),
\ \ W_n'(0;\gamma^2),
\ \ W_n''(0;\gamma^2)\ \  \mbox{and} 
\ \ W_n'''(0;\gamma^2),
\label{alg2:values}
\end{equation}
which give the initial conditions for (\ref{eqs:appell}).
Since most of  the solutions of Appell's equation
are highly oscillatory, and we are seeking a solution
which is not, it is necessary to use a solver which
is well-suited for ``stiff'' ordinary differential equations.
We use a fairly standard spectral method
 whose result is a  a  piecewise Chebyshev expansion of $W_n(x;\gamma^2)$ 
given on a partition of $[0,1)$.
We once again took the order of our expansion to be $k=29$,
and the partition on $[0,1)$  is determined through an adaptive procedure
reminiscent of the algorithm of Section~\ref{section:adap}.

The function $\PsiS{n}'(x;\gamma^2)$ is related to 
$W_n(x;\gamma^2)$ via
\begin{equation}
\PsiS{n}'(x;\gamma^2)  = \frac{1}{W_n(x;\gamma^2)},
\end{equation}
and we use this relation to construct a $k$th order
piecewise Chebyshev expansion of 
$\PsiS{n}'(x;\gamma^2)$ on the interval $[0,1)$.  The
value of  $\PsiS{n}(0;\gamma^2)$ is known --- in fact,
\begin{equation}
\PsiS{n}(0;\gamma^2) = - \frac{\pi}{2} (n+1),
\label{alg2:constant}
\end{equation}
and a $k$th order 
piecewise Chebyshev expansion of 
$\PsiS{n}(x;\gamma^2)$ on $[0,1)$
is obtained through the spectral integration of 
$\PsiS{n}'(x;\gamma^2)$ over $[0,1)$ with
(\ref{alg2:constant}) providing the constant
of integration.

Once the   $k$th order piecewise Chebyshev expansions of  $\PsiS{n}(x;\gamma^2)$
and $\PsiS{n}'(x;\gamma^2)$ are obtained, the 
function $\ps{n}(x;\gamma^2)$ can be evaluated at any point
$x$ in the interval $[0,1)$ by evaluating
these Chebyshev expansions at $x$ and then applying the formula
\begin{equation}
\ps{n}(x;\gamma^2) = 
\frac{\sin\left(\PsiS{n}(x;\gamma^2)\right)}{\sqrt{\PsiS{n}'(x;\gamma^2)}}\sqrt{1-x^2}.
\end{equation}
The function of second kind can also be evaluated, if it is so desired,
via
\begin{equation}
\qs{n}(x;\gamma^2) = 
\frac{\cos\left(\PsiS{n}(x;\gamma^2)\right)}{\sqrt{\PsiS{n}'(x;\gamma^2)}}\sqrt{1-x^2}.
\end{equation}

\end{section}

\begin{section}{Numerical experiments}
\label{section:experiments}

In this section, we describe numerical experiments conducted to evaluate
the performance of the algorithm of this paper.
Our code was written in Fortran and  compiled with the GNU
Fortran compiler version 7.4.0. 
Our implementation of the algorithm of this paper and 
our code for conducting the numerical  experiments described here is 
available on GitHub at the following
address:
\begin{center}
\url{https://github.com/JamesCBremerJr/Prolates}
\end{center}
All calculations were carried out  on an Intel Xeon E5-2697 processor 
running at 2.6 GHz.

In our implementation of the Xiao-Rokhlin algorithm, which is included 
the software mentioned above, the dimension of the tridiagonal
symmetric discretization matrix is taken to be
$n + \sqrt{n \gamma}$.
That is, we take the hidden constant in (\ref{introduction:cost})  to be $1$.
We found this to be sufficient to achieve near double precision accuracy.

\begin{subsection}{The Sturm-Liouville eigenvalues $\chi_n(\gamma^2)$}
\label{section:experiment1}

In these experiments, we measured the speed and 
accuracy with which our expansions 
 evaluate $\chi_n(\gamma^2)$
via comparison with the Xiao-Rokhlin algorithm.
In each experiment, $250,000$ pairs of the parameters
$\gamma$ and $n$ were constructed by 
choosing $500$ equispaced values of $\gamma$ in a specified
range and then, for each chosen value of $\gamma$,
picking $500$ random values of $n$ in 
the range $200 \leq n \leq  \gamma$.
For each pair of the parameters generated in this way, the eigenvalue  $\chi_n(\gamma^2)$ 
was evaluated via the expansion of Section~\ref{section:algorithm1}
and via the Xiao-Rokhlin algorithm.

Table~\ref{table:chi} reports the results of these experiments.  Each row there
corresponds to one experiment, and hence one range of values of $\gamma$.  
The values of   $\chi_n(\gamma^2)$ produced by the two algorithms
were compared at a total of  $3,000,000$ points
during the course of these experiments.

\begin{table}[h!!]
\small
\begin{center}
\begin{tabular}{cccc}
\toprule
Range of $\gamma$        &Maximum relative         &Average time             &Average time             \\
                         &difference               &expansion                &Xiao-Rokhlin             \\
\midrule
256 - 512 & 2.57\e{-15} & 8.23\e{-07} & 4.74\e{-04}  \\
512 - 1\sep,024 & 1.91\e{-15} & 7.80\e{-07} & 6.12\e{-04}  \\
1\sep,024 - 2\sep,048 & 2.07\e{-15} & 7.76\e{-07} & 1.14\e{-03}  \\
2\sep,048 - 4\sep,096 & 1.98\e{-15} & 6.79\e{-07} & 2.19\e{-03}  \\
4\sep,096 - 8\sep,192 & 2.09\e{-15} & 6.78\e{-07} & 4.36\e{-03}  \\
8\sep,192 - 16\sep,384 & 2.15\e{-15} & 6.94\e{-07} & 8.86\e{-03}  \\
16\sep,384 - 32\sep,768 & 1.64\e{-15} & 6.89\e{-07} & 1.79\e{-02}  \\
32\sep,768 - 65\sep,536 & 2.06\e{-15} & 6.85\e{-07} & 3.64\e{-02}  \\
65\sep,536 - 131\sep,072 & 2.21\e{-15} & 7.19\e{-07} & 7.48\e{-02}  \\
131\sep,072 - 262\sep,144 & 2.76\e{-15} & 7.14\e{-07} & 1.66\e{-01}  \\
262\sep,144 - 524\sep,288 & 4.93\e{-15} & 7.26\e{-07} & 3.71\e{-01}  \\
524\sep,288 - 1\sep,048\sep,576 & 6.40\e{-15} & 7.34\e{-07} & 1.05\e{+00}   \\
\bottomrule
\end{tabular}

\end{center}
\caption{A comparison of the time required to compute the 
Sturm-Liouville eigenvalue $\chi_n(\gamma)$ using
the method of this paper and via the Xiao-Rokhlin algorithm.
All times are in seconds.  Each row of the table corresponds
to $250,000$ evaluations of $\chi_n(\gamma)$.}
\label{table:chi}
\end{table}

\end{subsection}


\begin{subsection}{The functions $\ps{n}(x;\gamma^2)$}

In these experiments, we measured the speed and accuracy with which
the algorithm of this paper evaluates the functions $\ps{n}(x;\gamma^2)$
via comparison with the Xiao-Rokhlin algorithm.  
In each experiment, $10,000$ pairs of the parameters
$\gamma$ and $n$ were constructed by 
choosing $100$ equispaced values of $\gamma$ in a specified
range and then, for each chosen value of $\gamma$,
picking $100$ random values of $n$ in 
the range $200 \leq n \leq  \gamma$.
For each pair of the parameters generated in this way, the 
function $\ps{n}(x;\gamma^2)$ was evaluated at $1,000$ 
points  using the algorithm of this paper 
and via Xiao-Rokhlin method.  

\begin{table}[h!!]
\begin{center}
\small
\begin{tabular}{ccc}
\toprule
Range of $\gamma$                         &Average   &Average    \\
                             &precomp time                     &precomp time                    \\
                         &phase algorithm          &Xiao-Rokhlin             \\
\midrule
256 - 512 & 2.36\e{-04} & 4.38\e{-04}  \\
512 - 1\sep,024 & 2.43\e{-04} & 6.91\e{-04}  \\
1\sep,024 - 2\sep,048 & 2.73\e{-04} & 1.22\e{-03}  \\
2\sep,048 - 4\sep,096 & 2.97\e{-04} & 2.34\e{-03}  \\
4\sep,096 - 8\sep,192 & 3.20\e{-04} & 4.53\e{-03}  \\
8\sep,192 - 16\sep,384 & 3.36\e{-04} & 9.59\e{-03}  \\
16\sep,384 - 32\sep,768 & 3.59\e{-04} & 1.87\e{-02}  \\
32\sep,768 - 65\sep,536 & 3.84\e{-04} & 3.71\e{-02}  \\
65\sep,536 - 131\sep,072 & 3.96\e{-04} & 7.94\e{-02}  \\
131\sep,072 - 262\sep,144 & 4.24\e{-04} & 2.02\e{-01}  \\
262\sep,144 - 524\sep,288 & 4.34\e{-04} & 5.00\e{-01}  \\
524\sep,288 - 1\sep,048\sep,576 & 4.52\e{-04} & 1.09\e{+00}   \\
\bottomrule
\end{tabular}

\end{center}
\caption{A comparison of the average time taken by the precomputation
step of our algorithm with the average time take by the precomputation
step of the Xiao-Rokhlin algorithm.  All times are in seconds. 
Each row of the table corresponds to the construction of $10,000$
nonoscillatory phase functions/Legendre expansions.
}
\label{table:ps1}
\end{table}

Tables~\ref{table:ps1} and \ref{table:ps2} present
the results. 
Each row of these tables correspond to one experiment and hence
one range of $\gamma$.   
Table~\ref{table:ps1} gives the 
average time required to compute the phase function
$\PsiS{n}(x;\gamma^2)$ using the algorithm
of Section~\ref{section:algorithm2} and compares with it the average time required
by the precomputation phase of the Xiao-Rokhlin algorithm.
Table~\ref{table:ps2} compares the average time required evaluate $\ps{n}(x;\gamma^2)$
at a single point via the nonoscillatory phase function
$\PsiS{n}(x;\gamma^2)$ produced by the algorithm of this paper
and using the Legendre expansion produced by the Xiao-Rokhlin algorithm, 
as well as the maximum observed absolute error in the 
value produced by the algorithm of this paper.
The values of $\ps{n}(x;\gamma^2)$ produced by the two algorithms
were compared at a total of $120,000,000$ points.

\begin{table}[h!!]
\begin{center}
\small
\begin{tabular}{cccc}
\toprule
Range of $\gamma$        &Maximum absolute         &Average       &Average       \\
                         &error                    &evaluation time                     &evaluation time                    \\
                         &                         &phase algorithm            &Xiao-Rokhlin             \\
\midrule
256 - 512 & 9.04\e{-14} & 1.41\e{-07} & 1.29\e{-05}  \\
512 - 1\sep,024 & 1.37\e{-13} & 1.40\e{-07} & 2.23\e{-05}  \\
1\sep,024 - 2\sep,048 & 2.13\e{-13} & 1.39\e{-07} & 4.07\e{-05}  \\
2\sep,048 - 4\sep,096 & 3.02\e{-13} & 1.39\e{-07} & 7.78\e{-05}  \\
4\sep,096 - 8\sep,192 & 4.73\e{-13} & 1.39\e{-07} & 1.50\e{-04}  \\
8\sep,192 - 16\sep,384 & 6.74\e{-13} & 1.39\e{-07} & 2.99\e{-04}  \\
16\sep,384 - 32\sep,768 & 8.71\e{-13} & 1.39\e{-07} & 5.94\e{-04}  \\
32\sep,768 - 65\sep,536 & 1.32\e{-12} & 1.44\e{-07} & 1.20\e{-03}  \\
65\sep,536 - 131\sep,072 & 1.76\e{-12} & 1.39\e{-07} & 2.51\e{-03}  \\
131\sep,072 - 262\sep,144 & 3.62\e{-12} & 1.41\e{-07} & 5.56\e{-03}  \\
262\sep,144 - 524\sep,288 & 8.09\e{-12} & 1.39\e{-07} & 1.23\e{-02}  \\
524\sep,288 - 1\sep,048\sep,576 & 6.81\e{-12} & 1.39\e{-07} & 2.71\e{-02}  \\
\bottomrule
\end{tabular}

\end{center}
\caption{A comparison of the average time taken to evaluate
$\ps{n}(x;\gamma^2)$ via the algorithm of this paper
and by  via the Xiao-Rokhlin algorithm.  All times are in seconds. 
Each row of the table corresponds to $10,000,000$ evaluations
of $\ps{n}(x;\gamma^2)$.
}
\label{table:ps2}
\end{table}

\end{subsection}

\end{section}

\begin{section}{Acknowledgments}
The authors thank Vladimir Rokhlin for many useful conversations regarding
this work and for providing his code for evaluating prolate spheroidal wave functions.
Funding for this work was provided  by National Science Foundation grant DMS-1418723,
and by a UC Davis Chancellor's Fellowship.
\end{section}

\begin{section}{References}
\bibliographystyle{acm}
\bibliography{prolates}
\end{section}


\end{document}